\documentstyle[A4]{article}
\begin{document}
\newcommand{\IR}{\mbox{I \hspace{-0.2cm}R}}
\centerline{\Huge Twisted Classical Phase Space}
$$
$$
\begin{center}{\large P. Czerhoniak and A. Nowicki}\\
$$
$$
{\it Institute of Physics, Pedagogical University,\\
pl. S{\l}owia\'{n}ski 6, 65-029 Zielona G\'{o}ra, Poland}
\end{center}
$$
$$
\begin{abstract}
We consider relativistic phase space constructed by the twist
procedure from the translation sector of the standard, nondeforned
Poincar\'{e}        algebra. Using the concept of cross product
algebra we derive two kinds of phase space with noncommuting
configuration space. The generalized uncertainty relations are
formulated.
\end{abstract}
$$
$$
{\bf 1.}\quad In \cite{LNT} the twisted classical Poincar\'{e}
algebras were described using two-tensor procedure. Recently,
some efforts are made to describe deformed, noncommutative
phase spaces with noncommuting configuration
space, via the duality between quantum algebras and quantum
groups \cite{M}. This  construction of phase space leans
on the concept of a cross product algebra (Heisenberg double).

It appears that we can start from the standard, nondeformed
(classical) Poincar\'{e} algebra, making a twist and apply
duality to twisted Poincar\'{e} algebra in order to obtain
relativistic  phase space with noncommuting positions.\\
This point of view is presented in this paper
and we show that in the simplest case we obtain the
configuration noncommuting space with the structure of two
dimensional inhomogenous Euclidean or hiperbolic algebras.

At the begining we summarize the basic notions for  our
considerations.\medskip

$\bullet$\quad{\bf Twisted coproduct} Let ${\cal A}=(A,\Delta,
S,\varepsilon )$ be the Hopf algebra with comultiplication
$\Delta$, antipode $S$ and counit $\varepsilon$.\\
if there exists an invertible function $F=\sum_i f_i\otimes
f^i\in {\cal A}\otimes {\cal A}$ such that\medskip
$$
\Delta^F(a)=F\cdot \Delta (a)\cdot F^{-1}
\eqno(1.1)
$$
determines  new comultiplication, or equivalently:
$$
(\Delta^{F}\otimes 1 )\Delta^{F} = (1\otimes \Delta^{F})
\Delta^{F} \qquad \qquad coassociativity
\eqno(1.2)
$$
Comultiplication $\Delta^{F}$ is called the {\bf twisted coproduct}.
\medskip

$\bullet$\quad Two Hopf algebras ${\cal A}=(A,\Delta,S,\varepsilon )$
and ${\cal A}^F=(A,\Delta^F,S^F,\varepsilon )$ are related by {\bf
two-tensor twisting} if there exists an invertible function
$$
  F=\sum_i f_i\otimes f^i\in {\cal A}\otimes {\cal A}
\eqno(1.3a)
$$
satisfying
$$
F_{23}(1\otimes\Delta )F = F_{12}(\Delta\otimes 1)F\qquad
(\varepsilon \otimes 1)F = (1 \otimes \varepsilon )F = 1
\eqno(1.3b)
$$
where
$$
F_{12} = F\otimes 1 \qquad F_{23} = 1\otimes F \qquad F_{13} = \sum_i
f_i\otimes 1\otimes f^i
\eqno(1.3c)
$$

$\bullet$\quad For the complex simple Lie algebras $\hat{g}$ the
twist function $F$ is  given by
$$
F=\exp f \quad\quad\quad f\in \hat{c}\otimes \hat{c}
\eqno(1.4)
$$
where $\hat{c}$ is the commutative subalgebra of $\hat{g}$ {\it
(Cartan or Borel subalgebra)}.\bigskip\\
Let ${\cal{P}}$ be an algebra, ${\cal{X}}$ a vector space.

$\bullet$\quad A {\bf left action} {\it(representation)} of
${\cal{P}}$ on ${\cal{X}}$  is a linear map
$$
\triangleright : {\cal{P}}\otimes {\cal{X}}\rightarrow
{\cal{X}} : p\otimes x\rightarrow p\triangleright x
\eqno(1.5a)
$$
such that
$$
(p\tilde{p})\triangleright x = p\triangleright (\tilde{p}
\triangleright x)\quad 1\triangleright x = x
\eqno(1.5b)
$$

$\bullet$\quad We say that ${\cal{X}}$ is a {\bf left
${\cal{P}}$-module}.\\
In the case where ${\cal{X}}$ is an algebra and
${\cal{P}}$ a bialgebra ({\it or Hopf algebra})

$\bullet$\quad We say that ${\cal{X}}$ is a {\bf left
${\cal{P}}$-module algebra} if
$$
p\triangleright (x\tilde{x}) = (p_{(1)}\triangleright x) (p_{(2)}
\triangleright \tilde{y}) \quad p\triangleright 1=\epsilon(p) 1
\eqno(1.6)
$$
where $\epsilon$ denotes counit and we use the
Sweedler's notation
\begin{displaymath}
\Delta(p) = \sum p_{(1)}\otimes p_{(2)}
\end{displaymath}
Let ${\cal{P}}$ be a Hopf algebra and ${\cal{X}}$ a
left ${\cal{P}}$-module algebra.

$\bullet$\quad {\bf Left cross product algebra}({\it smash product})
${\cal{X}} >\!\!\!\triangleleft{\cal{P}}$  is a vector
space ${\cal{X}}\otimes {\cal{P}}$ with product \cite{M}
$$
(x\otimes p)(\tilde{x}\otimes \tilde{p})=x(p_{(1)}\triangleright
\tilde{x})\otimes p_{(2)}\tilde{p}\qquad
(left\,\,\, cross\,\,\, product)
\eqno(1.7)
$$
with unit element $1\otimes 1$, where $x, \tilde{x}\in
{\cal{X}}$ and $p, \tilde{p}\in {\cal{P}}$.\medskip

$\bullet$\quad {\bf commutation relations in a cross product algebra}.
\noindent The obvious isomorphism ${\cal{X}}\sim {\cal{X}}\otimes 1$,
${\cal{P}}\sim
1\otimes {\cal{P}}$ gives us the following cross relations between
the algebras ${\cal{X}}$ and ${\cal{P}}$
$$
[x,p]=x\circ p-p\circ x \qquad where \quad x\circ p=x\otimes p
\quad p\circ x=(p_{(1)}\triangleright x)\otimes p_{(2)}
\eqno(1.8)
$$
$$
$$

{\bf 2.}\quad We consider following two choices of the commuting
subalgebra $\hat{c}$ giving noncommuting configuration
space \cite{LNT}.\medskip\\
{\bf i)}\quad Let  $\hat{c}=(M_3=M_{12}, P_3, P_0)$, then
the twist function is given by $(r,s=3,0)$:
$$
F=e^{M_{3}\otimes A_{1}}e^{P_{r}\otimes B_{r}}=e^{A_{2}\otimes
M_{3}}e^{C_{r}\otimes P_{r}}
\eqno(2.1)
$$
where:\medskip\\
\hspace*{4em}$A_1 = \alpha_+ M_3 + (\delta_{+}^{r} + \delta_{-}^{r})
P_r \qquad B_r=(\delta_{+}^{r} - \delta_{-}^{r})M_3 +
\rho_{+}^{rs}P_s$\medskip\\
\hspace*{4em}$A_2= \alpha_+ M_3 + (\delta_{+}^{r} - \delta_{-}^{r})P_r
\qquad C_r=\rho_{+}^{sr}P_s + (\delta_{+}^{r} +
\delta_{-}^{r})M_3$\hfill$(2.2)$\medskip\\
Using the formula $(1.1)$ we obtain the twisted coproduct in the
form\medskip\\
$\Delta^F(M_\pm)=M_\pm \otimes e^{\pm A_1} + e^{\pm A_2}\otimes M_\pm
\pm P_\pm \otimes B_3e^{\pm A_1} \pm e^{\pm A_2}C_3 \otimes P_\pm$
\medskip\\ $\Delta^F(M_3)=M_3\otimes 1 + 1 \otimes M_3\hfill(2.3a)$
\bigskip\\$\Delta^F(N_\pm)= N_\pm \otimes e^{\pm A_1}
+ e^{\pm A_2} \otimes N_\pm -i{P_\pm \otimes B_0 e^{\pm A_1}
+ C_0 e^{\pm A_2}\otimes P_\pm}$\medskip\\
$\Delta^F(N_3)=N_3 \otimes 1 + 1 \otimes N_3 -i{P_3 \otimes B_0 + C_0
\otimes P_3 + P_0 \otimes B_3 + C_3 \otimes P_0}\hfill(2.3b)$\bigskip\\
$\Delta^F(P_1)=P_1 \otimes \cosh(A_1)+ \cosh(A_2) \otimes P_1
+ iP_2 \otimes \sinh(A_1) + i\sinh(A_2) \otimes P_2$\medskip\\
$\Delta^F(P_2)=P_2 \otimes \cosh(A_1) + \cosh(A_2) \otimes P_2
-iP_1 \otimes \sinh(A_1) -i\sinh(A_2)\otimes P_1$\medskip\\
$\Delta^F(P_3)=P_3 \otimes 1 + 1 \otimes P_3$\medskip\\
$\Delta^F(P_0)=P_0 \otimes 1 + 1 \otimes P_0\hfill(2.3c)$\bigskip\\
{\bf ii)}\quad For the second choice of the commuting subalgebra
$\hat{c}=(P_1, P_2, N_3=M_{30})$ the twist function takes the form
(a,b=1,2):
$$
F=e^{N_{3}\otimes A_{1}}e^{P_{a}\otimes B_{a}}=e^{A_{2}\otimes
N_{3}}e^{C_{a}\otimes P_{a}}
\eqno(2.4)
$$
where:\medskip\\
\hspace*{4em}$A_1 =(\xi_{+}^{a}-\xi_{-}^{a})P_{a}+\gamma_{+}N_{3}
\qquad B_a=(\rho_{+}^{ab}P_{b}+(\xi_{+}^{a}+\xi_{-}^{a})N_{3}$
\medskip\\ \hspace*{4em}$A_2= (\xi_{+}^{a}-\xi_{-}^{a})P_{a}
+\gamma_{+}N_{3} \qquad
C_a=\rho_{+}^{ab}P_{b}+(\xi_{+}^{a}+\xi_{-}^{a})N_{3}\hfill(2.5)$
\bigskip\\and the following formulae for coproduct
can be easy obtain:\medskip\\
$\Delta^F(M_\pm)=M_\pm \otimes \cos(A_1) + \cos(A_2)\otimes M_\pm
\pm \{N_\pm \otimes \sin(A_1)+\sin(A_2)\otimes N_\pm\}$\medskip\\
\hspace*{5em}$\mp\{P_3 \otimes (B_1\pm iB_2)\cos(A_1) + (C_1\pm
iC_2)\cos(A_2)\otimes P_3\}$\medskip\\
\hspace*{5em}$\mp i\{P_0\otimes (B_1\pm iB_2)\sin(A_1) + (C_1 \pm
iC_2)\sin(A_2)\otimes P_0\}$\medskip\\
$\Delta^F(M_3)=M_3 \otimes 1 + 1\otimes M_3$\medskip\\
\hspace*{5em}$ - i\{P_2
\otimes B_1 + C_1 \otimes P_2\} + i\{P_1\otimes B_2 + C_2\otimes
P_1\}\hfill(2.6a)$\bigskip\\
$\Delta^F(N_\pm)= N_\pm \otimes \cos(A_1) + \cos(A_2)\otimes N_\pm
\mp \{ M_\pm \otimes \sin(A_1) + \sin(A_2) \otimes M_\pm \}$\medskip\\
\hspace*{5em}$-i \{ P_0 \otimes (B_1 \pm iB_2)\cos(A_1) + (C_1 \pm
iC_2)\cos(A_2) \otimes P_0 \}$\medskip\\
\hspace*{5em}$ +P_3\otimes (B_1 \pm iB_2)\sin(A_1) + (C_1 \pm iC_2)
\sin(A_2) \otimes P_3$\medskip\\
$\Delta^F(N_3)= N_3 \otimes 1 + 1 \otimes N_3\hfill(2.6b)$\bigskip\\
$\Delta^F(P_\pm)= P_\pm \otimes 1 + 1 \otimes P_\pm$\medskip\\
$\Delta^F(P_3)= P_3 \otimes \cos(A_1) + \cos(A_2) \otimes P_3
+ iP_0 \otimes \sin(A_1) + i\sin(A_2) \otimes P_0$\medskip\\
$\Delta^F(P_0)=P_0 \otimes \cos(A_1) + \cos(A_2) \otimes P_0
+ iP_3 \otimes \sin(A_1) + i\sin(A_2) \otimes P_3\hfill(2.6c)$
\medskip\\ where:
$$
B_1 \pm iB_2 = (\rho_{+}^{1b} \pm i \rho_{+}^{2b})P_b
+(\xi_{+}^1 \pm i\xi_{+}^2 - (\xi_{-}^1 \pm i\xi_{-}^2))N_3
\eqno(2.7a)
$$
$$
C_1 \pm iC_2=(\rho_{+}^{1b} \pm i\rho_{+}^{2b})P_b
+ (\xi_{+}^1 \pm i\xi_{+}^2 + (\xi_{-}^1\pm i\xi_{-}^2))N_3
\eqno(2.7b)
$$
in both cases we assume the following hermiticity condition for
coproduct
$$
(\Delta^F)^{+} = \Delta^F \Rightarrow A^{+}_1=-A_1\,,\,
A^{+}_2=-A_2
\eqno(2.8)
$$
which leads to the hermitean configuration space.\bigskip

{\bf 3.}\quad Let us consider the first choice of the commuting
subalgebra (i). The hermicity condition $(8)$ give us
$$
A_1 = i(\alpha_+ M_3 + (\delta^{r}_{+} + \delta^{r}_{-})P_r\quad
A_2 = i(\alpha_{+} M_3 + (\delta^{r}_{+} - \delta^{r}_{-})P_{r})\quad
\alpha_{+}, \delta^{r}_{+},\delta^{r}_{-} \in \IR
\eqno(3.1)
$$
In the simple case $\alpha_{+}\equiv 0, \delta^{r}_{+}\equiv 0$ the
coproduct for the momentum takes the form\medskip\\
$\hspace*{2em}\Delta^F (P_r) = P_r \otimes I + I\otimes P_r$\medskip\\
$\hspace*{2em}\Delta^F(P_1) = P_1 \otimes
\cos(\tilde{A})+\cos(\tilde{A})\otimes P_1 + P_2 \otimes
\sin(\tilde{A})-\sin(\tilde{A})\otimes P_2$\medskip\\
$\hspace*{2em}\Delta^F(P_2) = P_2 \otimes \cos(\tilde{A})+
\cos(\tilde{A}) -P_1 \otimes \sin(\tilde{A}) +
\sin(\tilde{A})\otimes P_1\hfill(3.2)$\bigskip
$$
\tilde{A}=\tilde{A}_1=-\tilde{A}_2=\delta^{r}_{-}P_r
$$
We obtain the twisted classical phase space using $(3.2),
(1.8)$ and duality relations\medskip
$$
<x_\mu, P_\nu>=-i\hbar g_{\mu\nu}\qquad\qquad g_{\mu\nu}=(-1,1,1,1)
$$
and we get the nonvanishing commutator relations\medskip\\
- {\it noncommuting configuration space}\bigskip\\
\hspace*{2em}$[x_0,x_3]=0$\medskip\\
\hspace*{2em}$[x_0,x_1]=2i\hbar \delta^{0}_{-} x_2\qquad\,\,\,\,
[x_3,x_1]=-2i\hbar \delta^{3}_{-}x_2$\medskip\\
\hspace*{2em}$[x_0,x_2]=-2i\hbar \delta^{0}_{-} x_1\qquad [x_3,x_2]
=2i\hbar \delta^{3}_{-}x_1\hfill(3.3a)$\medskip\\
- {\it position-momentum cross relations}\bigskip\\
\hspace*{2em}$[p_1,x_0]=-i\hbar \delta^{0}_{-}p_2\qquad\quad
[p_2,x_0]=i\hbar \delta^{0}_{-} p_1$\medskip\\
\hspace*{2em}$[p_1,x_3]=i\hbar \delta^{3}_{-}p_2\qquad\quad\,\,\,\,
[p_2,x_3]=-i\hbar \delta^{3}_{-} p_1$\medskip\\
\hspace*{2em}$[p_1,x_1]=-i\hbar\cos(\tilde{A})\qquad
[p_2,x_1]=i\hbar\sin(\tilde{A})$\medskip\\
\hspace*{2em}$[p_1,x_2]=-i\hbar\sin(\tilde{A})\qquad
[p_2,x_2]=-i\hbar\cos(\tilde{A})\hfill(3.3b)$\medskip\\
and commuting momentum space.

We see that the twisting of nondeformed, standard Poincar\'{e}
algebra provides the noncommuting configuration space ({\it
noncommuting position operators}) for suitable choice of commuting
subalgebra $\hat{c}$. The simplest structure of noncommutative
configuration space $\{x_0,x_1,x_2\}$ is realized by two
dimensional inhomogenous Euclidean algebra
$iso(2)$ if we assume $\delta^{3}_{-}\equiv 0,
\alpha=\delta^{0}_{-}$.\\Then we get the commutation
relations in the form\medskip\\
\hspace*{2em}$[x_0,x_1]=2i\hbar \alpha x_2\qquad\qquad\,\,\,
[x_0,p_1]=i\hbar \alpha p_2$\medskip\\
\hspace*{2em}$[x_0,x_2]=2i\hbar \alpha x_1\qquad\qquad\,\,\,
[x_0,p_2]=-i\hbar \alpha p_1$\hfill(3.4a)\bigskip\\
\hspace*{2em}$[p_1,x_1]=-i\hbar\cos(\alpha p_0)\qquad
[p_1,x_2]=-i\hbar\sin(\alpha p_0)$\medskip\\
\hspace*{2em}$[p_2,x_1]=i\hbar\sin(\alpha p_0)\qquad\,\,\,\,\,
[p_2,x_2]=-i\hbar\cos(\alpha p_0)\hfill(3.4b)$\bigskip\\
other commutators vanish.

The phase space obtained in this way has a natural realization
by the standard quantum-mechanical momentum and position operators.
In fact, let $\hat{x}_\mu$ and $\hat{p}_\nu$ satisfy the Heisenberg
commutation relations $[\hat{x}_\mu,\hat{p}_\nu]=i\hbar
g_{\mu\nu}\qquad g_{\mu\nu}=(-1,1,1,1)$, then we define\medskip
$$
x_0 = \hat{x}_0 + \alpha(\hat{x}_1 \hat{p}_2 - \hat{p}_1\hat{x}_2)
= \hat{x}_0 +\alpha M_3
\eqno(3.5a)
$$
$\qquad\left( \begin{array}{c}
x_1\\ x_2
\end{array} \right) =\left( \begin{array}{cc}
\cos(\alpha\hat{p}_0) & -\sin(\alpha\hat{p}_0)\\
\sin(\alpha\hat{p}_0) & \cos(\alpha\hat{p}_0)
\end{array} \right)\,\,\,\left( \begin{array}{c}
\hat{x}_1\\ \hat{x}_2
\end{array} \right)\hfill(3.5b)$\bigskip\\
Therefore, the noncommuting configuration space one can obtain as
$\hat{p}_0$ - depending rotation in two dimensional commuting space
$(\hat{x}_1,\hat{x}_2)$.

Introducing the dispersion of observable
$a$ in quantum mechanical sense by\medskip
$$
\Delta(a) \ = \ \sqrt{<a^2> - <a>^2} \qquad
\Delta(a)\Delta(b)\geq {1\over 2}|<c>|
\eqno(3.6)
$$
where\,\,  c=[a,b], we obtain the generalized Heisenberg uncertainty
relations in the form\medskip\\
\hspace*{2em}$\Delta(x_0)\Delta(x_1)\geq\hbar \alpha |<x_2>|\qquad
\qquad \Delta(x_0)\Delta(p_1)\geq{1\over 2}\hbar \alpha |< p_2>|$
\medskip\\ \hspace*{2em}$\Delta(x_0)\Delta(x_2)\geq\hbar \alpha
|<x_1>|\qquad \qquad\Delta(x_0)\Delta(p_2)\geq{1\over 2}\hbar
\alpha |<p_1>|$\hfill(3.7a)\bigskip\\
\hspace*{2em}$\Delta(p_1)\Delta(x_1)\geq{1\over 2}\hbar |<\cos(\alpha
p_0)>|\qquad \Delta(p_1)\Delta(x_2)\geq{1\over 2}\hbar|<\sin(\alpha
p_0)>|$\medskip\\ \hspace*{2em}$\Delta(p_2)\Delta(x_1)\geq{1\over 2}
\hbar |<\sin(\alpha p_0)>|\qquad \Delta(p_2)\Delta(x_2)
\geq{1\over 2}\hbar|<\cos(\alpha
p_0)>|\hfill(3.7b)$\bigskip\\
Let us notice, that the standard Heisenberg relations we obtain
in the limit $\alpha\to 0$ or for the quantized energy
$\alpha p_0 = n\pi\quad n=0,\pm1,\pm2,...$.\bigskip

{\bf 4}\quad The second choice (ii) of the commuting subalgebra
$\hat{c}$ one can consider analogously to the case (i).
Using the hermicity condition $(2.8)$ and
the relations $(2.5)$ and $(2.6c)$, assuming for simplicity
$\gamma =0, \xi^{a}_{+}=0$ we obtain ($a,b=1,2$)\medskip\\
- {\it noncommuting configuration space}\bigskip\\
\hspace*{2em}$[x_0,x_a]=-2i\hbar \xi^{a}_{-}x_3$\medskip\\
\hspace*{2em}$[x_3,x_a]=-2i\hbar \xi^{a}_{-}x_0\hfill(4.1a)$\medskip\\
- {\it position-momentum cross relations}\bigskip\\
\hspace*{2em}$[p_0,x_0]=i\hbar\cosh(\xi^{a}_{-}p_a)
\qquad [p_3,p_0]=-i\hbar\sinh(\xi^{a}_{-} p_a)$\medskip\\
\hspace*{2em}$[p_0,x_1]=-i\hbar\xi^{1}_{-} p_3\qquad\qquad
[p_3,x_1]=-i\hbar\xi^{1}_{-} p_0$\medskip\\
\hspace*{2em}$[p_0,x_2]=-i\hbar\xi^{2}_{-} p_3\qquad\qquad
[p_3,x_2]=-i\hbar\xi^{2}_{-} p_0$\medskip\\
\hspace*{2em}$[p_0,x_3]=i\hbar\sinh(\xi^{a}_{-}p_a)\qquad
[p_3,x_3]=-i\hbar\cosh(\xi^{a}_{-}p_a)\hfill(4.1b)$\bigskip\\
and commuting momentum space.\medskip\\
Similary to the Euclidean case (i) the simplest structure of
noncommutative configuration space $\{x_0,x_1,x_2\}$ is
realized by two dimensional
inhomogenous hiperbolic algebra $iso(1,1)$ assuming
$\xi^{2}_{-}\equiv 0, \beta=\xi^{1}_{-}$. Now, we obtain the
commutation relations in the form\bigskip\\
\hspace*{2em}$[x_0,x_1]=-2i\hbar \beta x_3$\medskip\\
\hspace*{2em}$[x_3,x_1]=-2i\hbar\beta x_0\hfill(4.2a)$\bigskip\\
\hspace*{2em}$[p_0,x_1]=-i\hbar\beta p_3\qquad\qquad
[p_3,x_1]=-i\hbar\beta p_0$\medskip\\
\hspace*{2em}$[p_0,x_3]=i\hbar\sinh(\beta p_1)\qquad
[p_3,x_0]=-i\hbar\sinh(\beta p_1)$\medskip\\
\hspace*{2em}$[p_0,x_0]=i\hbar\cosh(\beta p_1)\qquad
[p_3,x_3]=-i\hbar\cosh(\beta p_1)\hfill(4.2b)$\bigskip\\
This kind of phase space one can also realized using
quantum-mechanical momentum and position operators
(see $(3.5)$)\medskip
$$
x_1 = \hat{x}_1 + \beta(\hat{x}_3 \hat{p}_0 - \hat{p}_3\hat{x}_0)
= \hat{x}_1 +\beta N_3
\eqno(4.5a)
$$
$\qquad\left( \begin{array}{c}
x_0\\ x_3
\end{array} \right) =\left( \begin{array}{cc}
\cosh(\beta\hat{p}_1) & \sinh(\beta\hat{p}_1)\\
\sinh(\beta\hat{p}_1) & \cosh(\beta\hat{p}_1)
\end{array} \right)\,\,\,\left( \begin{array}{c}
\hat{x}_0\\ \hat{x}_3
\end{array} \right)\hfill(4.5b)$\bigskip\\
therefore, noncommuting configuration space
one can obtain as $\hat{p}_1$ - depending
hiperbolic rotation in	two dimensional
commuting space $(\hat{x}_0,\hat{x}_3)$.\bigskip\\
The choice (ii) of twisting provides the generalized Heisenberg
uncertainty relations as follow\medskip\\
$\Delta(x_0)\Delta(x_1)\geq\hbar \beta|< x_3>|$\bigskip\\
$\Delta(x_3)\Delta(x_1)\geq\hbar\beta|< x_0>|\hfill(4.6a)
$\bigskip\\
$\Delta(p_0)\Delta(x_1)\geq{1\over 2}\hbar\beta|< p_3>|
\qquad\qquad \Delta(p_3)\Delta(x_1)\geq{1\over 2}\hbar\beta
|< p_0>|$\medskip\\ $\Delta(p_0)\Delta(x_3)\geq{1\over 2}\hbar|<
\sinh(\beta p_1)>|\qquad \Delta(p_3)\Delta(x_0)\geq{1\over 2}\hbar|
<\sinh(\beta p_1)>|$\medskip\\
$\Delta(p_0)\Delta(x_0)\geq{1\over 2}\hbar|<\cosh(\beta
p_1)>|\qquad \Delta(p_3)\Delta(x_3)\geq{1\over 2}\hbar|<\cosh(\beta
p_1)>|\hfill(4.6b)$
$$
$$
We see that two generalizations of the Heisenberg uncertainty
relations Eqs.$(3.7)$, $(4.6)$ have different behaviour for high
momentum limit.\\
In the case (i) they are bounded by ${\hbar\over 2}$ because of sine
and cosine functions. However in the hiperbolic case (ii) they are
strongly divergent for $p_1 \to\infty$.

\end{document}